%% file: HladkyPauls-CMC.tex
\documentclass{amsart}[11pt]
\usepackage{myAMSpackages}
\usepackage{Universal}
\usepackage{MyAMSEnvironments}
\usepackage{graphicx,subfigure,epsfig}

\newcommand{\wM}{\ensuremath{\widetilde{M}}\xspace}

\begin{document}
\title{Constant Mean curvature surfaces in sub-Riemannian geometry}
\author{Robert K. Hladky}
\address{Dartmouth College, Hanover, NH 03755}
\email{robert.hladky@dartmouth.edu}
\author{Scott D. Pauls}
\address{Dartmouth College, Hanover, NH 03755}
\email{scott.pauls@dartmouth.edu}
\thanks{Both authors are partially supported by NSF grant DMS-0306752}

\begin{abstract}
We investigate the minimal and isoperimetric surface problems in a
large class of sub-Riemannian manifolds, the so-called Vertically
Rigid spaces.  We construct an adapted connection for such spaces and,
using the variational tools of Bryant, Griffiths and Grossman, derive
succinct forms of the Euler-Lagrange equations for critical points for
the associated variational problems.  Using the Euler-Lagrange
equations, we show that minimal and isoperimetric surfaces satisfy a
constant horizontal mean curvature conditions away from characteristic
points.  Moreover, we use the formalism to construct a horizontal
second fundamental form, $II_0$, for vertically rigid spaces and, as a first
application, use $II_0$ to show that minimal surfaces cannot have
points of horizontal positive curvature and, that minimal surfaces in Carnot groups cannot be locally
horizontally geometrically convex.  We note that the convexity
condition is distinct from others currently in the literature.  
\end{abstract}

\maketitle
\input{intro}

\input{VerticallyRigid}
\input{BGG}

\input{ContactManifold}
\input{CMC}
\input{II}

\bibliographystyle{plain}
\bibliography{References}

\end{document}

%% file: intro.tex
\section{Introduction}\setS{intro}
Motivated by the classical problems of finding surfaces of least area
among those that share a fixed boundary (the minimal surface problem) and surfaces of least area
enclosing a fixed volume (the isoperimetric problem), several authors have recently formulated and
investigated similar problems in the setting of sub-Riemannian or Carnot-Carath\'edory spaces.  In
particular, N. Garofalo and D.M. Nheiu in \cite{GarNh} laid the
foundations of the theory of minimal surfaces in Carnot-Carath\'eodory
spaces and provided many of the variational tools necessary to make
sense of such a problem.  Building on this foundation, Danielli,
Garofalo and Nhieu, \cite{DGN}, investigated aspects of minimal and
constant mean curvature surfaces in Carnot groups.  Among many other results, these
authors showed the existence of isoperimetric sets and that, when considering the isoperimetric problem in the Heisenberg group, if one
restricts to a set of surfaces with a generalized cylindrical symmetry, then the
minimizers satisfy an analogue of the constant mean curvature
equation.  In this setting, the authors identify the absolute
minimizer bounding a fixed volume and show that it is precisely the
surface that Pansu conjectured to be the solution to the isoperimetric
problem (\cite{Pansu:isoper}).  Further in this direction of the
isoperimetric problem, Leonardi and
Rigot, \cite{LR}, independently showed the existence of isoperimetric sets
 in any Carnot group and investigated some of their properties.
 Leonardi and Masnou, \cite{LM}, investigated the geometry of the
 isoperimetric minimizers in the Heisenberg group and also showed a
 version of the result in \cite{DGN} showing the among sets with a
 cylindrical symmetry, Pansu's set is the isoperimetric minimizer.

In addition to this more general work in Carnot groups and
Carnot-Carath\'eodory spaces, a great deal of work has been done on
the minimal surface problem in more specialized settings.  For
example, the second author, in \cite{Pauls:minimal}, showed a connection
between Riemannian minimal graphs in the Heisenberg group and those in the
Carnot Heisenberg group and used this connection to prove $W^{1,p}$ estimates for solutions to the minimal surface equation.  In
addition, he found a number of initial examples of minimal surfaces in
the Heisenberg group and used them to demonstrate non-uniqueness of the
solution to the Dirichlet problem in the Heisenberg group.  
Recently, both Garofalo and the second author ,\cite{GP}, and Cheng, Huang,
Malchiodi and Yang, \cite{CHMY}, independently investigated minimal
surfaces in some limited settings.  Garofalo and the second author restricted
their view to the Heisenberg group and provided, among other results,
a representation theorem for smooth minimal surfaces, a horizontal regularity theorem and proved an
analogue of the Bernstein theorem, showing that a minimal surface in
the Heisenberg group that is a graph over some plane satisfies a type
of constant curvature condition.  We note that Cheng and Huang, 
\cite{CH}, independently showed a more general version of this
Bernstein-style theorem by classifying all properly embedded minimal
surfaces in $\mathbb{H}$.  In \cite{CHMY}, the authors investigate
$C^2$ minimal surfaces in three dimensional pseudohermitian geometries
(including, of course, the Heisenberg group) and, using the techniques
of CR geometry, investigate the structure of minimal surfaces in this
setting and, among many results, prove, under suitable conditions, a
uniqueness result for the Dirichlet problem for minimal surfaces in
the Heisenberg group.  In \cite{Pauls:Obstr}, the second author extends the
representation theorem of \cite{GP} to $C^1$ minimal surfaces in
$\mathbb{H}^1$, provides examples of continuous (but not smooth) minimal surfaces and
shows a geometric obstruction to the existence of smooth minimal
solutions to the Plateau Problem in the Heisenberg group.  In
\cite{Cole}, Cole examines minimal surfaces in spaces of
Martinet-type.  While this collection of spaces includes the
Heisenberg group and many of those considered in \cite{CHMY}, Cole's
thesis also treats three dimensional spaces that do not have
equiregular horizontal subbundles.  In \cite{Cole}, Cole derives the
minimal surface equation and examines the geometry and existence of
smooth solutions to the Plateau Problem.  

While there has been great progress in the understanding of minimal
and constant mean curvature surfaces in the setting of
Carnot-Carath\'eodory spaces, there are still many fundamental open
questions left to address.  Most notably, much of the focus has been
on the minimal surface equation and the majority of the work has
focused on more limited settings such as the Heisenberg group, groups
of Heisenberg type or three dimensional pseudohermitian
manifolds.  In this paper, we will address more general problems using a new tool to help discriminate between
the various types of constant mean curvature surfaces that abound in
different Carnot-Carath\'eodory spaces. 

\vspace{.25in}
\noindent
{\bf Question:}  {\em In a Carnot-Carath\'eodory manifold $M$, do the surfaces of
  least perimeter or the surfaces of
  least perimeter enclosing a fixed volume satisfy any partial
  differential equations?  Can the solutions be characterized geometrically?}

\vspace{.25in}

In Euclidean space, there is a beautiful connection between the
geometry of surfaces and the solutions to these variational problems:
minimal surfaces are characterized as zero mean curvature surfaces
while isoperimetric surfaces have constant mean curvature.   

In this paper, we restrict ourselves to a large class of sub-Riemannian
manifolds which we call vertically rigid sub-Riemannian (VR) spaces.  These
spaces are defined in \rfS{VR} and include basically all examples already studied (including Carnot groups,
Martinet-type spaces and pseudohermitian manifolds) but is a much larger class.  On such spaces, we define a new connection,
motivated by the Webster-Tanaka connection of strictly pseudoconvex pseudohermitian
manifolds, that is adapted to the sub-Riemannian structure.  Using
this connection and the variational framework of  Bryant, Griffiths
and Grossman (\cite{BGG1}), we investigate minimal and isoperimetric surface
problems.  The framework of \cite{BGG1} provides a particularly nice
form of the Euler-Lagrange equations for these problems and leads us to define a {\em horizontal second fundamental form},
$II_0$, and the {\em horizontal mean curvature}, Trace($II_0$), associated to a hypersurface in a Carnot group.  Given a
noncharacteristic submanifold $\Sigma$ of a VR space $M$, let $\{e_1,\dots, e_k\}$ be an
orthonormal basis for for the horizontal portion of the tangent space
to $\Sigma$ and let $e_0$ be the unit horizontal normal to $\Sigma$ (see the
next section for precise definitions).  Then, we define the horizontal
second fundamental form as 
\[II_0 = \begin{pmatrix}\langle\nabla_{e_1}e_0,e_1\rangle & \dots &
    \langle\nabla_{e_1}e_0,e_k\rangle\\
\vdots & \vdots & \vdots \\
\langle\nabla_{e_k}e_0,e_1\rangle& \dots &\langle\nabla_{e_k}e_0,e_k\rangle\end{pmatrix}\]
and define the horizontal mean curvature as the trace of
$II_0$.  We note that the notion of horizontal mean curvature
has appeared in several contexts (see \cite{BC}, \cite{CHMY}, \cite{DGN},
 \cite{DGN:convex}, \cite{GP}, \cite{Pauls:minimal}), and that this notion
coincides with the others, possibly up to a constant multiple.
However, we emphasize that the version of the mean curvature above
applies to all VR spaces (before this work, only \cite{DGN} deals with mean
curvature in any generality, but again is limited to Carnot groups) and has the advantage of being written 
in an invariant way with respect to the fixed surface.  With
this notion in place, we have a characterization of $C^2$ solutions to the
two variational problems discussed above:

\bgT{Main 1} Let $M$ be a vertically rigid sub-Riemannian manifold and
  $\Sigma$ a noncharacteristic $C^2$ hypersurface.  $\Sigma$ is  a
  critical point of the first variation of perimeter if and only if
  the horizontal mean curvature of $S$ vanishes.
\enT

Similarly,

\bgT{Main 2} Let $M$ be a vertically rigid sub-Riemannian space and
  $\Sigma$ a $C^2$ hypersurface.
  If $\Sigma$ is a  solution to the isoperimetric problem,  then the horizontal mean curvature of $\Sigma$ is locally constant.
\enT
   
Thus, we recover an analogue of the classical situation:  the
solutions to these two problems are found among the critical points of
the associated variational problems.  Moreover, these critical points
are characterized by having the trace of the second fundamental form
be constant. 

We note that the characterization of minimal surfaces in terms of both
a PDE and in terms of mean curvature was
achieved first by Danielli, Garofalo and Nhieu, \cite{DGN}, in Carnot
groups and by the second author, \cite{Pauls:minimal}, for minimal graphs in
the Heisenberg group.  The technique described above provides a broad
extension of this characterization and describes mean curvature in a
geometrically motivated manner.  From this point of view, this is most
similar to the treatment of mean curvature by Cheng, Huang, Malchiodi and Yang,
\cite{CHMY}, who use a similar formalism but used the Webster-Tanaka
connection.  In contrast, some of the 
earlier definitions of mean curvature relied on the minimal surfaces
equation for the definition (as in \cite{Pauls:minimal}) or via a
different geometric analogue such as a symmetrized horizontal Hessian
(as in \cite{DGN}, \cite{DGN:convex}).   

On the other hand, for
isoperimetric surfaces, the only known links between isoperimetric
sets and constant mean curvature were under the restriction of
cylindrical symmetry 
in the first Heisenberg group (\cite{DGN},\cite{LM}), for $C^2$
surfaces in 3-dimensional pseudohermitian manifolds (\cite{CHMY}) and
for $C^2$ surfaces in the 
first Heisenberg group using mean curvature flow methods due to Bonk and
Capogna (\cite{BC}) and, recently, for $C^2$ surfaces in all
Heisenberg groups due to Rigor\'e and Rosales (\cite{RR}).  Thus, our
treatment of these problems unifies these results and extends them to
a much larger class of sub-Riemannian manifolds.  Moreover, we provide
a number of new, general techniques and tools for investigating these
problems in a very general setting.  

As mentioned above,
the geometric structure of minimal surfaces has only been studied in
cases such as the Heisenberg group, pseudohermitian manifolds, and
Martinet-type spaces.  In general, even in the higher Heisenberg
groups, nothing is known about the structure and geometry of minimal
surfaces.  As an illustration of the power of this framework, we use the
horizontal second fundamental form to provide some geometric
information about minimal surfaces in any VR space.  To better describe
minimal surfaces, we introduce some new notions of curvature in VR
spaces:

\bgD{hpc} Let $II_0$ be the horizontal second fundamental form for a
$C^2$ noncharacteristic surface, $\Sigma$, in a vertically rigid sub-Riemannian
manifold $M$. Let $\{\mu_0, \dots, \mu_k\}$ be the eigenvalues
(perhaps complex and with multiplicity) of $II_0$.  Then, the {\bf
  horizontal principle curvatures} are given by 
\[ \kappa_i = \real{\mu_i}\]
for $0 \le i \le k$.  

Moreover, given $x \in \Sigma$, we say that $\Sigma$ is {\bf horizontally 
positively (non-negatively) curved} at $x$ if $II_0$ is either positive (semi-)definite or negative (semi-) definite at $x$ and
is {\bf horizontally negatively curved} at $x$ if there is at least one positive and one negative $\kappa_i$.  $\Sigma$ is {\bf horizontally flat} at $x$ if $\kappa_i
=0$ for $0 \le i\le k$.  
\enD

Let $\Sigma$ be a $C^2$ hypersurface in $M$, a
vertically rigid sub-Riemannian manifold.  Then, the {\bf horizontal
exponential surface} at $x \in \Sigma$, $\Sigma_0(x)$, is defined to be the union of
all the horizontal curves in $\Sigma$ passing through $x$.  The notion
of horizontal principle curvatures described above gives rise to a new
definition of convexity:

\bgD{horiz convex}  A subset $U$ of a Carnot group $M$ with $C^2$ boundary $\Sigma$ is {\bf horizontally geometrically
  convex} (or hg-convex) if, at each noncharacteristic point $x \in \Sigma$,
$\Sigma_0(x)$ lies to one side of $T_x^h\Sigma$, the horizontal
tangent plane to $\Sigma$ at $x$.  We say that $\Sigma$ is locally
hg-convex at $x$ if there exists an $\epsilon >0$ so that $\Sigma_0(x)
\cap B(x,\epsilon)$ lies to one side of $T_x^h\Sigma$.
\enD

We note that this notion of convexity is distinct from those described
in \cite{DGN:convex} or \cite{LMS}.  In \rfS{II}, we give
explicit examples showing the nonequivalence of the various notions.

With these definitions in place, we prove an analogue to the classical
statement that a minimal surfaces in $\mathbb{R}^3$ must be nonpositively curved.

\bgT{Main 3} Let $\Sigma$ be a $C^2$ noncharacteristic minimal
hypersurface in a vertically rigid sub-Riemannian space $M$.  Then,
$\Sigma$ cannot contain a point of horizontal positive curvature.  If we further assume that $M$ is a Carnot group, then
$\Sigma$ cannot be locally horizontally geometrically convex. 
\enT

We emphasize that this is the first
description of the geometry of minimal surfaces in a relatively general class of spaces.  

The rest of the paper is divided into five sections.  In
\rfS{VR}, we define vertically rigid sub-Riemannian spaces, the
adapted connection we mentioned above, and an adapted frame bundle for
such objects.  In \rfS{BGG}, we briefly review the relevant machinery from \cite{BGG1}.  In \rfS{CM}, we address the minimal surface
problem using the machinery of Bryant, Griffiths and Grossman.
\rfS{CMC} addresses the isoperimetric problem in this setting
and finally, in \rfS{II}, we define the horizontal second
fundamental form and prove the geometric properties of minimal
surfaces described above.

%% file: VerticallyRigid.tex
\section{Vertically rigid sub-Riemannian manifolds}\setS{VR} 
We begin with our basic definitions:
\bgD{SR space}  A {\bf sub-Riemannian (or Carnot-Carath\'eodory) manifold}
is a triple $(M,V_0,\langle \cdot, \cdot \rangle)$ consisting of a
smooth manifold $M^{n+1}$, a smooth $k+1$-dimensional distribution $V_0
\subset TM$ and a smooth inner product on $V_0$.  This structure is
endowed with a metric structure given by
\[d_{cc} (x,y) = \inf \{\int \langle \dot{\gamma},\dot{\gamma}
\rangle^\frac{1}{2} | \gamma(0)=x,\gamma(1)=y, \gamma \in
\mathscr{A}\}\]
where $\mathscr{A}$ is the space of all absolutely continuous paths
whose derivatives, when they are defined, lie in $V_0$.
\enD

\bgD{Vertically Rigid} A sub-Riemannian manifold has a {\bf vertically
  rigid complement} if there exist
\begin{itemize}
\item a smooth complement $V$ to $V_0$ in $TM$
\item a smooth frame $T_1,\dots T_{n-k}$ for $V$
\item a Riemannian metric $g$ such that $V$ and $V_0$ are orthogonal, $g$ agrees with $\langle \cdot, \cdot \rangle$ on $V_0$ and $T_1,\dots T_{n-k}$ are orthonormal.
\item a partition of $\{1..,n-k\}$ into equivalence classes such that for all sections $X \in \Gamma(V_0)$, $g ( [X,T_j],T_i)=0$ if $j \sim i$.
\end{itemize}

A sub-Riemannian space with a vertically rigid complement is called a
{\bf vertically rigid (VR) space}.
\enD

For a VR space, we shall denote the number of equivalence classes of the partition by $v$ and the size of the partitions by $l_1,\dots l_v$. In particular, we then have $l_1+\dots+l_v=n-k$. After choosing an order for the partitions, for $j >0$ we set
\[ V_j = \text{span} \{T_i: \text{ $i$ is in the $j$th partition.}\}\]
Then \[ TM = \bigoplus_{j=0}^v V_j.\]
After reordering we can always assume that the vector fields $T_1,\dots T_{l_1}$ span $V_1$, the next $l_2$ span $V_2$ etc.

There are 3 motivating examples for this definition:
\bgX{Pseudo}
Let $(M,\theta,J)$ be a strictly pseudoconvex pseudohermitian structure (see \cite{Tanaka}). Then $V_0=\ker{\theta}$ has codimension $1$ and a vertically rigid structure can be defined by letting $T_1$ be the characteristic (Reeb) vector field of $\theta$ and defining $g$ to be the Levi metric
\[ g(X,Y) = d\theta(X,JY) + \theta(X)\theta(Y).\]
Since $T_1$ is dual to $\theta$ and $T_1 \lrcorner d\theta =0$ the required commutation property clearly holds.

\enX

\bgX{Carnot}
Let $(M,\mathfrak{v}_0)$ be a graded Carnot group with step size $r$. Then the Lie algebra of left-invariant vector fields of $M$ decomposes as
\[ \mathfrak{m} = \bigoplus_{j=0}^r \mathfrak{v}_j \]
where $\mathfrak{v}_{j+1}=[\mathfrak{v}_0,\mathfrak{v}_j]$ for $j<r$ and $[\mathfrak{v}_0,\mathfrak{v}_r]=0$. We then set $V_j = \text{span} (\mathfrak{v}_j)$ and  construct a (global) frame of left invariant vector fields for each $\mathfrak{v}_j$, $j>0$. This clearly yields a vertically rigid structure.

\enX

In the previous examples, the vertical structure was chosen to carefully mimic the bracket-generating properties of the sub-Riemannian distribution. We include another example, where the bracket-generation step size need not be constant to illustrate the flexibility of this definition.

\bgX{Martinet}
Let $M=\rn{3}$. We define a Martinet-type sub-Riemannian structure on $M$ by defining $V_0$ to be the span of
\[ X= \partial_x + f(x,y) \partial_z, \qquad Y=\partial_y+g(x,y) \partial_z\] where $f$ and $g$ are smooth functions. The metric is defined by declaring $X,Y$ an orthonormal frame for $V_0$. (In particular if we take $f=0$ and $g=x^2$, we see that the step size is $1$ on $x \ne 0$ and $2$ at $x=0$.) Now define $T_1 = \partial_z$ and extend the metric so that $X,Y,T_1$ are orthonormal. Again the commutation condition clearly holds.

\enX

The illustrate the generality of the definition, we give one last example.

\bgX{vfds}  Let $\{X_1,\dots,X_k\}$ be a collection of smooth vector fields
on $\mathbb{R}^n$ that satisfy H\"ormander's condition.  We will construct the
$\{T_i\}$ as follows.  As the $X_i$ bracket generate, let $\{T_i\}$ be
a basis for the complement of the span of the $\{X_i\}$ formed by
differences of the brackets of the $X_i$ and linear combinations of
the $X_i$ themselves.  These $T_i$ are naturally graded by counting
the number of brackets of $X_i's$ it takes to include $T_i$ in the
span.  Define a Riemannian inner product that makes the
$\{X_1,\dots,X_k,T_1,\dots, T_{n-k}\}$ an orthonormal basis.  This
structure satisfies all the conditions for a vertically rigid
structure except possibly the last.  

We note that the majority of the examples in the literature, either
from subelliptic PDE, control theory and/or robotic path planning, 
satisfy the last condition.  

\enX

The advantage of vertically rigid structures is that they admit connections which are adapted to analysis in the purely horizontal directions. 

\bgD{Adapted}
A connection $\nabla$ on $TM$ is adapted to a vertically rigid structure if
\begin{itemize}
\item $\nabla$ is compatible with $g$, i.e. $\nabla g=0$.
\item $\nabla T_j =0$ for all $j$.
\item $\text{Tor}_p(X,Y) \in V_p$ for all sections $X,Y$ of $V_0$ and $p \in M$.
\end{itemize} 
\enD 

The motivating example for this definition is the Webster-Tanaka connection for a strictly pseudoconvex pseudohermitian manifold \cite{Tanaka}.

\bgL{Existence}
Every vertically rigid structure admits an adapted connection.
\enL

\pf
Let $\overline{\nabla}$ denote the Levi-Cevita connection for $g$. Define $\nabla$ as follows: set $\nabla T_j =0$ for all $j$. Then for a section $X$ of $V_0$ and any vector field $Z$ define
\[ \nabla_Z X = (\overline{\nabla}_Z X)_0\]
where $(\cdot)_0$ denotes the orthogonal projection onto $V_0$. This
essentially defines all the Christoffel symbols for the connection. It
is easy to see that it satisfies all the required conditions.  For example, to show $\nabla g=0$, take vector fields $X,Y,Z$ and write
\begin{equation*}
\begin{split}
X &= X_0 + \sum x_i T_i\\
Y &= Y_0 + \sum y_i T_i \\
\end{split}
\end{equation*}
where $X_0,Y_0 \in V_0$.  Using the fact that $\nabla T_i=0$, we have
\begin{equation*}
\begin{split}
\nabla g (X,Y,Z) &= Zg(X,Y)-g(\nabla_ZX,Y)-g(X,\nabla_ZY)\\
&=Zg(X_0,Y_0)+\sum Z(x_iy_i) -
g((\overline{\nabla}_ZX_0)_0,Y_0)\\&-g(\sum
Zx_i,Y)-g(X_0,(\overline{\nabla}_ZY_0)_0)-g(\sum Zy_i,X)\\
&=g(\overline{\nabla}_ZX_0,Y_0)-g((\overline{\nabla}_ZX_0)_0,Y_0)+g(X_0,\overline{\nabla}_ZY_0)-g(X_0,(\overline{\nabla}_ZY_0)_0)\\&+\sum
Z(x_iy_i)-\sum Z(x_i)y_i - \sum Z(y_i)x_i\\
&=0
\end{split}
\end{equation*}
The last equality follows since $X_0$ and $Y_0$ are horizontal vector
fields and using the product rule.   The statement about torsion
follows directly from the definition.

\epf

\bgL{unique}
If $X,Y,Z$ are horizontal vector fields then
\begin{align*}
 \langle \nabla_X Y, Z \rangle &= \frac{1}{2} \big( X\ap{Y}{Z}{} + Y\ap{X}{Z}{}-Z\ap{X}{Y}{}\\
 & \quad +\ap{[X,Y]_0}{Z}{}+\ap{[Z,X]_0}{Y}{}+\ap{[Z,Y]_0}{X}{} \big).\end{align*}
In particular, this depends solely on the choice of orthogonal complement $V$.\enL

\pf We note that since $\nabla g=0$  and $V_0$ is parallel,
\[ \ap{\nabla_X Y}{Z}{} = X\ap{Y}{Z}{} - \ap{Y}{\nabla_X Z}{}.\]
Since the torsion of two horizontal vector fields is purely vertical,
we also obtain 
\begin{align*}
 \ap{\nabla_X Y}{Z}{}&= \ap{\nabla_Y X}{Z}{} + g([X,Y],Z) \\
 &=\ap{\nabla_Y X}{Z}{} + \ap{[X,Y]_0}{Z}{}.\end{align*}
The remainder of the proof is identical to the standard treatment of the Levi-Cevita connection on a Riemannian manifold given in \cite{Chavel}.

\epf

\bgR{floppy}  We note that the definition of adapted connection leaves
some flexibility in its definition.  In particular, we have some
freedom in defining Christoffel symbols related to $\nabla_{T_i} X$ when $X$ is a
section of $V_0$.  While we could make choices that would fix a unique
adapted connection, we will not do so in order to preserve the maximum
flexibility for applications.
\enR

To study these connections and sub-Riemannian geometry it is useful to introduce the idea of the graded frame bundle.

\bgD{GradedFrame}
An orthonormal frame $(e,t)=e_0,\dots,e_k,t_1,\dots,t_{n-k}$ is graded if $e_0,\dots,e_k$ span $V_0$, $t_1,\dots,t_{n-k}$ span $V$ and each $t_j$ is in the span of $\{ T_i: i\sim j\}$.
\enD

The bundle of graded orthonormal frames $\begin{CD}\mathcal{GF}(M) @>\pi >> M \end{CD}$ is then a $O(k+1) \times \Pi_{j=1}^{v} O(l_j)$-principle bundle.

On $\mathcal{GF}(M)$ we can introduce the canonical $1$-forms $\omega^j$, $\eta^j$ defined at a point $f=(p,e,t)$ by
\[ \omega^j (X)_f = g_p( \pi_* X, e_j), \quad \eta^j (X)_f = g_p( \pi_* X, t_j).\] An adapted connection can be viewed as a affine connection on $\mathcal{GF}(M)$. The structure equations are then determined by the following lemma.

\bgL{Structure}
On $\mathcal{GF}(M)$ there exist connection $1$-forms $\omega_j^i$ $0 \leq i,j \leq k$ and $\eta^i_j$, $1 \leq i, j \leq n-k$ together with torsion $2$-forms $\tau^i$ $0 \leq i \leq k$ and $\tilde{\tau}^i$ $1 \leq i \leq n-k$ such that
\bgE{Structure}
\begin{split}
d\omega^i &= \sum\limits_{0 \leq j \leq k}\omega^j \wedge \omega_j^i + \tau^i,\\
d\eta^i &= \sum\limits_{i \sim j} \eta^j \wedge \eta_j^i + \tilde{\tau}^i.
\end{split}
\enE
\enL

\pf The content of the lemma is in the terms that do not show up from the standard structure equations of an affine connection. However since $V_0$ is parallel we can immediately deduce that there exist forms $\omega^j_k$ such that $\nabla e_k = \omega^j_k \otimes e_j$. Furthermore since each $t_j$ is in the span of $\{T_i: i \sim j\}$  and all the $T_i$ are also parallel we must have $\nabla t_j = \sum_{i \sim j} \eta^i_j \otimes t_i$ for some collections of forms $\eta^i_j$.  

\epf

\bgL{Properties} The torsion forms for an adapted connection have the following properties:
\begin{itemize}
\item $\tau^j(e_a,e_b)=0$
\item $\tilde{\tau}^j(t_i,e_b)=0$ if $j \sim i$
\end{itemize}
for any lifts of the vector fields.
\enL

\pf
The first of these is a direct rewrite of the defining torsion condition for an adapted connection. For the second we observe that 
\[ \text{Tor}(T_i,e_b)= \nabla_{T_i} e_b - \nabla_{e_b} T_i - [T_i,e_b]\]
is orthogonal to $T_j$ if $i \sim j$ by the bracket conditions of a vertically rigid structure. The result then follows from noting that torsion is tensorial.
  
\epf

%% file: BGG.tex
\section{Exterior differential systems and variational problems}\setS{BGG}

In this section, we briefly review the basic elements of the formalism
of Bryant, Griffiths and Grossman  which can be found in
more detail in chapter one of \cite{BGG1}.  Their formalism requires
the following data:

\begin{enumerate}
\item A contact manifold $(M,\theta)$ of dimension $2n+1$
\item An $n$-form, called the Lagrangian, $\Lambda$ and the associated
  area functional
\[ \mathscr{F}_\Lambda(N) = \int_N \Lambda\]
where $N$ is a smooth compact Legendre submanifold of $M$, possibly
with boundary.  In this setting, a Legendre manifold is a manifold
$i:N\rightarrow  M$ so that $i^* \theta =0$.
\end{enumerate}

From this data, we compute the Poincar\'e-Cartan form $\Pi$ from the
form $d\Lambda$.  They show that $d\Lambda$ can be locally expressed
as 
\[ d \Lambda = \theta \wedge (\alpha + d\beta) + d(\theta \wedge
\beta)\]
for appropriate forms $\alpha,\beta$.  Then,
\[\Pi = \theta \wedge (\alpha + d \beta)\]
and we often denote $\alpha + d\beta$ by $\Psi$.  With this setup,
Bryant, Griffiths and Gross prove the following characterization of
Euler-Lagrange systems (\cite{BGG1} section 1.2):
\bgT{BGG}  Let $N$ be a Legendre surface with boundary $\partial N$ in
$M$ given by $i:N \rightarrow  M$ as above.  Then $N$ is a
stationary point under all fixed boundary variations, measured with
respect to  $\mathscr{F}_\Lambda$, if and only if $i^*\Psi =0$.  
\enT
In the next two sections, we will use this formalism and the previous
theorem to investigate the minimal and isoperimetric surface
problems.

%% file: ContactManifold.tex
\section{Minimal Surfaces}\setS{CM}

For a  $C^2$ hypersurface $\Sigma$ of a vertically rigid sub-Riemannian manifold we define the horizontal perimeter of $\Sigma$ to be
\bgE{PerimeterG}
P(\Sigma) = \int\limits_\Sigma \left| (\nu_g)_0 \right| \nu_g \lrcorner dV_g
\enE
where $\nu_g$ is the unit normal to $\Sigma$ with respect to the Riemannian metric $g$. At noncharacteristic points of $\Sigma$, i.e. where $T\Sigma \nsubseteq V_0$, this can be re-written as
\bgE{Perimeter}
P(\Sigma)= \int\limits_\Sigma \nu \lrcorner dV_g
\enE
where $\nu$ is the horizontal unit normal vector, i.e. the projection of the Riemannian normal to $V_0$.  We note that, when
restricting to the class of $C^2$ submanifolds, this definition is
equivalent to the perimeter measure of De Giorgi introduced in
\cite{DeG}.   Our primary goal for this section is to answer the following question.

\bgQ{Minimal} 
In a vertically rigid sub-Riemannian manifold, given a fixed boundary can the hypersurfaces spanning the boundary with least perimeter measure be geometrically characterized?
\enQ 
   
To answer this question, we shall employ the formal techniques of Bryant, Griffiths and Grossman \cite{BGG1} by exhibiting the minimizing hypersurfaces as integrable Legendre submanifolds of a contact covering manifold $\wM$. More specifically, we define $\wM$ to be the bundle of horizontally normalized contact elements,
\[ \wM = \{ (p,\nu,T) \in M \times (V_0)_p \times V_p:  \|\nu\|=1\}.\]
Thus $\begin{CD} \wM @>\pi_1 >> M\end{CD}$ has the structure of an $\sn{k+1} \times \rn{n-k}$-bundle over $M$.  Next we define a contact form $\theta$ on $\wM$ by 
\[ \theta_{\tilde{p}}(X) = g_p( (\pi_1)_* X, \nu+T).\]
To compute with $\theta$ it is useful to work on the graded frame
bundle. However as there is no normalization on the $T$ component of
$\wM$, we shall need to augment $\mathcal{GF}(M)$ to the fiber bundle
$\mathcal{GF}^0(M)$ defined as follows: over each point the fibre is
$O(k+1) \times \Pi_{j=1}^{v} O(l_j) \times \rn{l_1+\dots + l_v}$. The
left group action is extend as follows.  If $h=(h_1,h_2) \in O(k+1)
\times  \Pi_{j=1}^{v} O(l_j)$, 
\[ h \cdot (p,e,t,a) = (p, (h_1 \cdot e,h_2 \cdot t), ah_2^{-1}).\]
The natural projection $\pi$ from $\mathcal{GF}^0(M)$ to $M$ now filters through $\wM$ as 
\[ \begin{CD} \mathcal{GF}^0(M) @>\pi_2>> \wM @>\pi_1>> M \end{CD}\]
where under $\pi_2$, $(p,e,t,a) \mapsto (p, e_0, \sum a_jt_j)$. In particular this means $\pi_2 \circ (id,h_2) = \pi_2$. This formulation now allows us to pull $\theta$ back to $\mathcal{GF}^0(M)$ by
\[ \pi_2^* \theta = \omega^0 + a_j \eta^j.\] We shall denote this pullback by $\theta$ also.

\bgR{Augment}
The augmented frame bundle $\mathcal{GF}^0(M)$ is not a principle bundle and so we cannot impose an affine connection on it in the usual sense. However since it has the smooth structure of $\mathcal{GF}(M) \times \rn{l_1+\dots +l_v}$ we can naturally include the canonical forms and the connection structure equations of $\mathcal{GF}(M)$ into the augmented bundle. Thus the results of \rfL[VR]{Structure} and \rfL[VR]{Properties} hold on $\mathcal{GF}^0(M)$ also.

\enR

\bgL{NonDegenerate}
The contact manifold $(\wM,\theta)$ is maximally non-degenerate, i.e. $\theta \wedge d\theta^{n}$.
\enL

\pf
We shall work on the augmented graded frame bundle where
\[ d\theta = \omega^j \wedge \omega_j^0 + \tau^0 + da_j \wedge \eta^j + a_j( \eta^i \wedge \eta_i^j + \tilde{\tau}^j).\]
We pick out one particular term of the expansion of $\theta \wedge d\theta^n$, namely  
\[ \mu = \omega^0 \wedge \omega^1 \wedge \dots \wedge \omega^k \wedge \omega^1_0 \wedge \dots \omega^k_0 \wedge \eta^1 \wedge \dots \eta^{n-k} \wedge da_1 \wedge \dots da_{n-k}.\]
The connection forms are vertical (in the principle bundle sense) and the canonical forms are horizontal (again in the bundle sense). Thus $\mu$ is the wedge of $n-k$ $da$ terms, $n+k+1$ horizontal forms and $k$ vertical forms. Since each torsion form is purely (bundle) horizontal, $\mu$ is clearly the only term of this form in $\theta \wedge d\theta^n$. All the forms are independent so $\mu$ does not vanish. Thus we deduce that $\theta \wedge d\theta^n \ne 0$ on $\mathcal{GF}^0(M)$ and so cannot vanish on $\wM$.

\epf

The transverse Legendre submanifolds of $(\wM,\theta)$ are the immersion $\iota\colon \Sigma \hookrightarrow \wM$ such that 
$\iota^*\theta =0$ and $\pi_2 \circ \iota$ is also an immersion. These are noncharacteristic oriented hypersurface patches in $M$ with normal directions defined by the contact element in $\wM$.

Define
\bgE{Lambda} \Lambda = \omega^1 \wedge \dots \wedge \omega^k \wedge \eta^1 \wedge \eta^{n-k}\enE on $\mathcal{GF}^0(M)$. Then $\Lambda = \pi_2^* ( \nu \lrcorner \pi_1^* dV)$ and so $\Lambda$ is basic over $\wM$. Furthermore due to \rfE{Perimeter} we see that
\bgE{PerimeterL} P(\Sigma)= \int_\Sigma \iota^* \Lambda.
\enE

Now on $\mathcal{GF}^0(M)$, \rfL[VR]{Properties} implies that $\tau^j$ has no component of the form $\omega^0 \wedge \omega^j$ and $\tilde{\tau}^j$ none of form $\omega^0 \wedge \eta^j$. Thus we see from \rfE[VR]{Structure} that
\begin{align*} d\Lambda&= \sum\limits_j  (-1)^{j-1} \omega^1 \wedge \dots \wedge (\omega^0 \wedge \omega^j_0) \wedge  \dots \wedge \omega^k \wedge \eta^1 \wedge \dots \wedge \eta^{n-k}\\
&  + \sum\limits_j  (-1)^{j-1}\omega^1 \wedge \dots \wedge (\omega^j \wedge \omega^j_j)  \wedge \dots \wedge \omega^k \wedge \eta^1 \wedge \dots \wedge \eta^{n-k}\\
&  +\sum\limits_j (-1)^{k+j-1}\omega^1 \wedge \dots \omega^k \wedge \eta^1 \wedge \dots \wedge (\eta^j \wedge \eta^j_j) \wedge \dots \wedge \eta^{n-k}.
\end{align*}
The connection is metric compatible so $\omega^j_j=0$ and $\eta^j_j=0$. Thus the second and third sums vanish identically. This implies $d\Lambda =\theta \wedge \Psi$ where
\bgE{Psi} \Psi= \sum_j \omega^1 \wedge \dots \wedge  (\omega^j_0) \wedge  \dots \wedge \omega^k \wedge \eta^1 \wedge \dots \wedge \eta^{n-k}.\enE
If $\Sigma \subset \wM$ is a transverse Legendre submanifold, then we can construct a graded frame adapted to $\Sigma$, i.e. with $e_0 =\nu$. Choosing any section immersing $\Sigma$ into $\mathcal{GF}^0(M)$ we can then pull $\Psi$ back to $\Sigma$. Switching the $\omega$'s and $\eta$'s to represent the coframe and connection form for this fixed frame, we get
\bgE{PsiSigma} \Psi_{|\Sigma} = \left(\sum_{j=1}^k \omega^j_0(e_j) \right) \; \omega^1 \wedge \dots \wedge \omega^k \wedge \eta^1 \wedge \dots \eta^{n-k}.\enE

\bgT{MinimalSurface}
Suppose $\Sigma$ is a $C^2$ hypersurface in the vertically rigid manifold $M$. Then $\Sigma$ is a critical point for perimeter measure in a noncharacteristic neighborhood $U \subset \Sigma$ if and only if the unit horizontal normal $\nu$ satisfies the minimal surface equation
\bgE{MSE2} H=0 \enE
in $U$. Equivalently, the horizontal normal must satisfy 
\bgE{MSE} \text{div}_g \nu =0\enE
everywhere on $U$, where the divergence is taken with respect to the Riemannian metric $g$. 
\enT

\pf
From the Bryant-Griffiths-Grossman formalism \cite{BGG1} we see that $\iota\colon\Sigma \hookrightarrow \wM$ is a stationary  Legendre submanifold for $\Lambda$ in a in a small neighborhood if and only if $\iota^* \Psi =0$. This condition is just $\sum_{j=1}^k \omega^j_0(\nu) =0$ for any local orthonormal frame $(\nu,e_1,\dots e_k)$ for $V_0$. This can be re-written as
\bgE{One} H=\sum \langle \nabla_{e_j} \nu ,e_j \rangle =0.\enE
A standard formula in Riemannian geometry (see for example \cite{Kobayashi}) states that for any connection for which the volume form is parallel, the divergence of a vector field can be computed by
\[ \text{div}_g X =  \text{trace}( \nabla X + \text{Tor}(X,\cdot) ).\] 
The adapted connection is symmetric for $g$ and so we can apply this result while noting that by the defining conditions $ {trace} (\text{Tor}(\nu,\cdot))=0$. Thus 
\[ \text{div}_g \nu = \sum \langle \nabla_{e_j} \nu, e_j \rangle + \sum  g( \nabla_{t_j} \nu, t_j) = \sum \langle \nabla_{e_j} \nu, e_j \rangle.\]
Referring back to \rfE{One} then completes the proof.

\epf

\bgC{Independence}
The minimal surface equation \rfE{MSE} may depend on the choice of orthogonal complement $V$, but not on the remainder of the vertically rigid structure or choice of adapted connection.\enC

\pf
After we write $\text{div}_g\nu = \sum \ap{\nabla_{e_j} \nu}{e_j}{}$, the result follows immediately from \rfL[VR]{unique}.
\epf

\bgC{Ruling}
Any minimal noncharacteristic patch of a vertically rigid sub-Riemannian manifold $(M,V_0,\langle \cdot, \cdot \rangle)$ with
\[ \text{dim }V_0 =2\]
 is ruled by horizontal $\nabla$-geodesics. 
\enC

\pf
Extend $\nu$ off $\Sigma$ to any unit horizontal vector field.  Define $\nu^\bot$ to be any horizontal unit vector field that is orthogonal to $\nu$. By the torsion properties of the connection and the arguments of \rfT{MinimalSurface} the minimal surface equation \rfE{MSE} can  be written
\[ 0 = \langle \nabla_{\nu^\bot} \nu, \nu^\bot \rangle = -  \langle \nu, \nabla_{\nu^\bot} \nu^\bot\rangle.\]
Since $\nu^\bot$ has no covariant derivatives in vertical directions, this implies that
\[ \nabla_{\nu^\bot} \nu^\bot =0.\]
In other words the integral curves of $\nu^\bot$ are $\nabla$-geodesics. However, these integral curves clearly foliate the noncharacteristic surface patch.

\epf

\bgR{generalization}  We note that the last corollary is a
generalization of the results of Garofalo and the second author
\cite{GP} in the Heisenberg group, those of Cheng, Huang, Malchiodi
and Yang \cite{CHMY} in three dimensional pseudohermitian manifolds
and those of Cole \cite{Cole} in Martinet-type spaces.
In those cases, the authors proved the minimal surfaces in those
settings were ruled by appropriate families of horizontal curves.  
\enR

%% file: CMC.tex
 \section{CMC surfaces and the isoperimetric problem}\setS{CMC}
 
 We now investigate the following question
 
 \bgQ{Isoperimetric}
 Given a fixed volume, what are the closed surfaces bounding this volume of minimal perimeter.
 \enQ
 
Using the results of the previous section, we can now define a
hypersurface of locally constant mean curvature (CMC) by requiring
that $H=\text{constant}$ on each connected component of $\Sigma\upp =
\Sigma - \text{char}(\Sigma)$.  If we wish to specify the constant, we
will call $\Sigma$ a CMC($\rho$) surface.  By comparing to the
Riemannian case, these are our prime candidates for solutions to
\rfQ{Isoperimetric}. Throughout this section we shall make the
standing assumption that the volume form $dV_g$ is globally exact, i.e
there exists a form $\mu$ such that $d \mu =dV_g$ on $M$. Since this
is always locally true, the results of this section will hold for
sufficiently small domains.   
 
For a closed codimension $2$ surface $\gamma$ in $M$ we define 
\bgE{Span}
\text{Span}(\gamma,a) = \left\{ C^2 \text{ noncharacteristic surface } \Sigma: \partial \Sigma = \gamma, \int_\Sigma \mu =a\right\}. 
\enE

\bgL{SpanMin}
If $\text{Span}(\gamma,a)$ is non-empty then any element of minimal perimeter $\Sigma_0$ must have constant mean curvature.
\enL

\pf
As the boundary of the surfaces are fixed, we can again employ the formalism of \cite{BGG1}. We permit variations that alter the integral $\int_\Sigma \mu$ and apply a Lagrange multiplier method to establish a condition for critical points of perimeter. Indeed the minimiser $\Sigma_0$ must be a critical point of the functional
\[ \Sigma \mapsto \int_\Sigma \Lambda - c \int_\Sigma \mu\] for some constant $c$.
Now pulled-back to the contact manifold
\begin{align*}
d (\Lambda -c \mu) &= \theta \wedge \Psi - c \pi_1^* dV_g\\
& = \theta  \wedge \left( (H -c)  \omega^1 \wedge \dots \wedge \omega^k \wedge \eta^1 \wedge \dots \wedge \eta^{n-k}\right)\\
& := \theta \wedge \widetilde{\Psi}. 
\end{align*}
The same methods as \rfT[CM]{MinimalSurface} then imply that $\widetilde{\Psi}_{|\Sigma_0}=0$ and hence $\Sigma_0$ has constant mean curvature.

\epf
 
 \bgT{CMC}
 If a $C^2$ domain $\Omega$ minimises surface perimeter over all domains with the same volume, then $\Sigma= \partial \Omega$ is locally CMC.
 \enT
 
 \pf
 Let $p \in \partial \Sigma$ be a noncharacteristic point. As $\Sigma$ is $C^2$ there exists a noncharacteristic neighbourhood $U$ of $p$ in $\Sigma$ with at least $C^2$ boundary. Now by Stokes' theorem $\text{Vol}(\Omega) = \int_\Sigma \mu$ and so $U$ must minimise perimeter over all noncharacteristic surfaces in $\text{Span}(\partial U, \int_U \mu)$. Therefore $U$ has constant mean curvature by \rfL{SpanMin}.
  
 \epf

Under certain geometric conditions we can provide a more intuitive description of $\mu$.
 \bgD{Dilation} A dilating flow for a vertically rigid structure is a global flow $F\colon M \times \rn{} \longrightarrow M$ 
 \begin{itemize}
  \item $(F_\lambda)_* E_j = e^{\lambda} E_j$ for some fixed horizontal orthonormal frame $E_1,\dots,E_{k+1}$.
  \item $(F_\lambda)_* T_j = e^{\gamma_j \lambda} T_j$ for some constant $\gamma_j$.   
 \end{itemize}
 Associated to a dilating flow are the dilation operators defined by
 \[ \delta_\lambda = F_{\log{\lambda}}\]
 and the generating vector field $X$ defined by
 \[ X_p = \frac{d}{d\lambda}_{|\lambda=0} F(\lambda,p).\]
 The homogeneous dimension of $M$ is given by \[ Q = k+1 + \sum\limits_{j=1}^{n-k}\gamma_j.\]
 The dilating flow is said to have an origin $O$ if for all $p$,  $\delta_\lambda (p) \to O$ as $\lambda  \to 0$.
 \enD
 
 In the sequel, a vertically rigid sub-Riemannian manifold that admits a dilating flow with origin will be referred to as a VRD-manifold.
 
 \bgX{Carnot}
 All Carnot groups admit a dilation with origin. On the Lie algebra
 level, the dilation is defined merely by defining a linear map with
 eigenspaces the various levels of the grading, i.e. $\delta_\lambda X
 = \lambda^{j+1} X$ for $X \in V_j$. The group dilations are then
 constructed by exponentiating.  
 \enX
 
 \bgX{Martinet}
 The jointly homogeneous Martinet spaces, i.e. those of
 \rfX[VR]{Martinet} with the functions $f$ and $g$ bihomogeneous of
 degree $m$. Then the dilations are defined by 
 \[ \delta_\lambda (x,y,z) = (\lambda x ,\lambda y, \lambda^{m+1} z).\]
 Then clearly \[(\delta_{\lambda})_* (\partial_x + f(x,y) \partial_z)
 =  \lambda \partial_x + \lambda^{m+1} f(x,y) \partial_z = \lambda(
 \partial_x + f(\lambda x, \lambda y) \partial_z)\] and
 $(\delta_{\lambda})_* \partial_z = \lambda^{m+1} \partial_z$.  
 \enX
 
 \bgL{Exact}
 In a VRD-manifold the form $ \mu =  Q^{-1} X \lrcorner dV_g$ satisfies \[d\mu =  dV_g.\]
 \enL
 
 \pf
 Let $\omega^j$, $\eta^i$ be the dual basis to $E_j$, $T_i$. Then
 \[ dV_g = \omega^1 \wedge \dots \wedge \omega^{k+1} \wedge \eta^1 \wedge \dots \wedge \eta^{n-k}.\] 
 Now $F_\lambda^* \omega^j(Y) = \omega^j((F_\lambda)_*Y)$ so $F_\lambda^* \omega^j = \lambda \omega^j$. Thus \[\mathcal{L}_X \omega^j = \frac{d}{d \lambda}_{|\lambda =0} F_\lambda^* \omega^j = \lambda \omega^j.\]
 By a virtually identical argument we see that $\mathcal{L}_X \eta^j = \gamma_j \eta^j$. Therefore 
 \[ d\mu = Q^{-1} d(X \lrcorner dV_g) = Q^{-1} \mathcal{L}_X dV_g =  dV_g.\]
  \epf

 Every point $p$ in a VRD-manifold can be connected to the origin by a curve of type $ t \mapsto \delta_t(p)$. For any surface $\Sigma$ we can then construct the dilation cone over $\Sigma$ as 
 \bgE{cc} \text{cone}(\Sigma) = \{ \delta_t(p): 0 \leq t \leq 1, p \in \Sigma.\}\enE
 
 \bgL{Volume}
 Suppose $\Sigma$ is $C^2$ surface patch in a VRD-manifold such that any dilation curve intersects $\Sigma$ at most once. If $\Sigma$ is oriented so that the normal points away from the origin then
 \[ \text{Vol}\left(\text{cone}(\Sigma)\right) =  \int\limits_\Sigma \mu.\]
 \enL
 
 \pf
 This is just Stokes' theorem for a manifold with corners, for
 \[ \int_{\text{cone}(\Sigma)} dV_g =  \int_{\partial \text{cone}(\Sigma)} \mu = \int_{\Sigma}  \mu\]
 as $\mu$ vanishes when restricted to any surface foliated by dilation curves.
  
 \epf
 
We can now interpret \rfL{SpanMin} as minimising surface perimeter under the constraint of fixed dilation cone volume. Since the volume of a domain is equal to the signed volume of its boundary dilation cone, this yields some geometric intuition for the arguments of \rfT{CMC}.

\bgR{RefDGN-LM}   In \cite{DGN} and \cite{LM}, the authors characterize cylindrically symmetric
  minimizers in the Heisenberg group as constant mean curvature
  surfaces in that setting.  Our treatment allows for such a
  characterization in all VR spaces without the assumption of
  cylindrical symmetry.  We note, however, that this method requires
  some regularity (the surfaces must be at least $C^2$ to ensure the
  computations work) while the work in \cite{DGN} and \cite{LM} is 
  more general in this respect, allowing for piecewise $C^1$ defining functions. \enR
\bgR{Rem-isop}  We point out that  \rfT{CMC} provides an approach
  to understanding the isoperimetric problem in VR or VRD spaces via a
  better understanding of their constant mean curvature surfaces.  For
  a specific example, the reader is referred to \rfS{II}
  below.
\enR

%% file: II.tex
\section{The horizontal second fundamental form}\setS{II}

We now present a more classical interpretation of these results by defining an analogue to the second fundamental form. 
\bgD{secfund} Consider a
noncharacteristic point of a hypersurface $\Sigma$ of a VR space $M$ and fix a
horizontal orthonormal frame $e_0,\dots e_k$ as before with $\nu
=e_0$. Then the {\bf horizontal second fundamental form} as the $k
\times k$ matrix, 
\bgE{h2}
II_0 = \begin{pmatrix} \ap{ \nabla_{e_1} \nu }{e_1}{} & \dots &   \ap{ \nabla_{e_1} \nu }{e_k}{}\\
\vdots & \vdots & \vdots \\
 \ap{ \nabla_{e_k} \nu }{e_1}{} & \dots &  \ap{ \nabla_{e_k} \nu }{e_k}{} \end{pmatrix} = \begin{pmatrix} \omega_0^1(e_1) &\dots& \omega_0^k(e_1)\\ \vdots & \vdots & \vdots \\
  \omega_0^1(e_k)& \dots &\omega_0^k(e_k)
 \end{pmatrix} .  
\enE

Further, we define the {\bf horizontal mean curvature}, $H$, to be the trace of
$II_0$.  
\enD

By arguments already
given we see $\text{div}_g\nu =H$. Thus we have shown that a $C^2$
surface that minimizes perimeter must satisfy $H=0$ at any
noncharacteristic point. 

We note that several authors have proposed other candidates for
certain types of analogues of the second fundamental form and
horizontal mean curvature.  For example, Danielli, Garofalo and Nhieu (\cite{DGN},\cite{DGN:convex}) use a symmetrized
horizontal Hessian to analyze minimal and CMC surfaces (in \cite{DGN})
and convex sets (in \cite{DGN:convex}).  We emphasize that many of the other candidates are symmetrized versions of the second fundamental form while the definition above is explicitly {\em a priori} non-symmetric. 

This definition of mean curvature coincides (up to a constant multiple) with the various
definitions of mean curvature in the Carnot group setting (see, for
example, \cite{DGN}, \cite{GP}).  Moreover, our version of the minimal
surface equation as $Trace \; II_0=0$ matches with others in the
literature once it is suitably interpreted.  
In particular, our formulation, when restricted to the appropriate setting is equivalent to that of \cite{DGN}, \cite{CHMY}, \cite{Pauls:minimal},
\cite{Cole}, and \cite{RR}.  

In addition to the horizontal mean curvature described above, we would
also like to
define other aspects of horizontal curvature.  

\bgD{hpc} Let $II_0$ be the horizontal second fundamental form for a
$C^2$ noncharacteristic surface, $\Sigma$, in a vertically rigid sub-Riemannian
manifold $M$. Let $\{\mu_0, \dots, \mu_k\}$ be the eigenvalues
(perhaps complex and with multiplicity) of $II_0$.  Then, the {\bf
  horizontal principle curvatures} are given by 
\[ \kappa_i = \real{\mu_i}\]
for $0 \le i \le k$.  

Moreover, given $x \in \Sigma$, we say that $\Sigma$ is {\bf horizontally 
positively (non-negatively) curved} at $x$ if $II_0$ is either positive (semi-)definite or negative (semi-) definite at $x$ and
is {\bf horizontally negatively curved} at $x$ if there is at least one positive and one negative $\kappa_i$.  $\Sigma$ is {\bf horizontally flat} at $x$ if $\kappa_i
=0$ for $0 \le i\le k$.  
\enD

This definition coupled with the observation that $Trace \; II_0 =
\kappa_0 + \dots + \kappa_k$ yields the following immediate corollary
of \rfT[CM]{MinimalSurface}:

\bgC{Negative Curvature}  If $\Sigma$ is a $C^2$  minimal surface in a
vertically rigid sub-Riemannian manifold $M$, then $\Sigma$ has no noncharacteristic 
points of horizontal positive curvature.
\enC

This is reflective of the Euclidean and Riemannian cases where
minimal surfaces cannot have points of positive curvature.

\bgR{definitelyrequired}  We note that having, for example, only
positive horizontal principle curvatures at a point is necessary but
not  sufficient to conclude
that the surface is horizontally positively curved.
\enR

\bgD{hcurve}  Let $c$ be a horizontal curve on $\Sigma$ a $C^2$ hypersurface
in a vertically rigid sub-Riemannian manifold.  Then at a noncharacteristic point of $\Sigma$, the {\bf horizontal curvature} of $c$ is given by
\[ k_c = <\nabla_{\dot{c}}\dot{c},e_0>\]

\enD

We note that, analogous to the Euclidean and Riemannian cases, there is a
connection between the horizontal curvature of curves passing through
a point on a hypersurface and the
horizontal second fundamental form at that point:

\bgL{hcurve connection}   Let $c$ be a horizontal curve on $\Sigma$ a $C^2$ hypersurface
in a vertically rigid sub-Riemannian manifold.  Then at noncharacteristic points,
\[k_c  = <II_0(\dot{c}),\dot{c}>\]
\enL
\pf Since $c$ is horizontal, we have that $\dot{c} = c_1\; e_1 + \dots + c_k
\; e_k$ for appropriate functions $c_i$.  Differentiating $<\dot{c},e_0> = 0$, we have:
\begin{equation*}
\begin{split}
<\nabla_{\dot{c}}\dot{c},e_0> &= - <\dot{c},\nabla_{\dot{c}} e_0>\\
k_c &= -<II_0(\dot{c}),\dot{c}>
\end{split}
\end{equation*}

\epf

This gives, as an immediate corollary, an analogue of Meusnier's theorem:
\bgC{Meusnier} All horizontal curves lying on a surface $\Sigma$ in $M$, a VR space, that, at a noncharacteristic point $x \in \Sigma$, have the same tangent vector also have the same horizontal curvature at $x$.
\enC
As in the classical case, \rfC{Meusnier} allows us to speak of the
horizontal curvature associated with a direction rather than with a
curve, showing that $II_0$ contains all of the horizontal curvature
information at a point.  

\bgL{Principle directions}  Given $\Sigma$ and $x$ as above, let $l$ be the number of distinct principle curvatures at $x$.  Then, there exist curves $\{c_1,\dots,c_l\}
\subset \Sigma$ so that 
\[k_{c_i} = \kappa_i\]
\enL
\pf  Let $\{\lambda_1, \dots \lambda_j,\lambda_{j+1}\pm i\beta_{j+1},
, \dots, \lambda_l\pm i
\beta_l\}$ be the eigenvalues of $II_0$ at $x$ associated with the distinct principle curvatures.  Further, let  $\{u_1,\dots,u_j,u_{j+1}\pm i v_{j+1},\dots,
u_k\pm i v_k\}$ be the associated eigenvectors.  Without loss of generality, we have ordered the
eigenvalues so that the real eigenvalues appear first and the complex
eigenvalues are last.  We note that, for each complex conjugate pair
of eigenvalues, $\lambda_j \pm \beta_j$, the associated principle
curvatures, $\kappa_j$ and $\kappa_{j+1}$, are equal.  Using the
eigenvectors, we may replace $\{e_1,\dots,e_k\}$
by a new basis given by
\[\{u_1,\dots,u_j,u_{j+1},v_{j+1},\dots,u_{l},v_{l},
\tilde{e}_{2l-j},\tilde{e}_k\}\] where $\{\tilde{e}_i\}$ form a basis
for the orthogonal complement of the eigenvectors.  Rewriting $II_0$ with respect to this
new basis, there is a submatrix of $II_0$ which is a block matrix
where the 
first block is of the form 
\[A =\begin{pmatrix}
\lambda_1& \dots & 0\\
\vdots &  \vdots & \vdots \\
0 & \dots & \lambda_j \\
\end{pmatrix}\]
and there are $k-j$ remaining blocks of the form
\[B_i=\begin{pmatrix}
\lambda_i & \beta_i\\
-\beta_i & \lambda_i
\end{pmatrix}\]

Now, let $c_i$ be the integral curve of the $i^\text{th}$ new basis
vector for $1<i<2l-j$.  Then,
\[k_{c_i} = -<II_0(\dot{c_i}),\dot{c_i}>= <II_0(e_i),e_i>=\kappa_i\]

\epf
\bgR{jordanform}  We note that since $II_0$ is often nonsymmetric,
there is often not a full basis of eigenvectors.  For example, in the
Carnot group $\mathbb{H} \times \mathbb{R}$ with coordinates
$(x,y,t,s)$ and Lie algebra spanned by $\{X_1,X_2,X_3,T\}$ where
\begin{equation*}
\begin{split}
X_1 &= \partial_x - \frac{y}{2} \partial_t\\
X_2 &= \partial_y + \frac{x}{2} \partial_t\\
X_3 &= \partial_s\\
X_4 &= \partial_t
\end{split}
\end{equation*}
Taking $V_0=span\{X_1,X_2,X_3\}$ and $V_1 = span \{X_4\}$ yields a
vertically rigid structure.  The surface defined by
$\frac{xy}{2}-t-s=0$ has unit horizontal normal given by:
\[ \nu = \frac{y}{\sqrt{1+y^2}} \; X_1 - \frac{1}{\sqrt{1+y^2}} \;
X_3\]
and 
\[II_0=\begin{pmatrix} 0 & \frac{y^2}{(y^2+1)^\frac{3}{2}} \\ 0 & 0 \end{pmatrix}\]
Thus, this is a minimal surface and $II_0$ has a double eigenvalue of $0$
and a single eigenvector $(1,0)$ in this basis.  The presents an
entirely different phenomena then the analogous Riemannian or
Euclidean situation.

We note that this phenomena and the existence of complex eigenvalues
both indicate the existence of a nontrivial bracket structure among
the elements of the tangent space to $\Sigma$.  Indeed, both of these
indicate that there are vector fields $e_i, e_j \in \{e_1,\dots,e_k\}$
with the property that $[e_i,e_j]$ has a component in the $e_0$
direction. In particular, this is an indication that the distribution
$\{e_1,\dots, \e_k\}$ is not integrable and hence, $\Sigma$ cannot be
realized as a surface ruled by codimension one horizontal submanifolds.
\enR
We pause to note that we can now state an analogue of
\rfC[CM]{Ruling}:

\bgC{CMC Ruling}  Any CMC($\rho$) noncharacteristic patch of a vertically rigid sub-Riemannian manifold $(M,H,\langle \cdot, \cdot \rangle)$ with
\[ \text{dim }V_0 =2\]
 is ruled by horizontal curves with horizontal constant curvature $\rho$. 
\enC
\pf
We follow precisely the same proof as that of \rfC[CM]{Ruling} to get
that 
\[  <\nabla_{\nu^\bot} \nu^\bot,\nu> =\rho\]
Thus, the integral curves of $\nu^\bot$ have constant horizontal
curvature $\rho$.  
\epf

\bgR{ExplicitRules}  We note that in specific cases, these rulings can
be computed exactly.  Basically this amounts to explicitly computing
the torsion terms and solving ODEs.  For example, it is
straightforward to verify that in the Heisenberg group, such curves
are geodesics with respect to the Carnot-Carath\'eodory metric and are
horizontal lifts of planar circles of curvature $\rho$.  This
particular observation is also contained in \cite{CHMY}.
\enR

\bgD{Horizontal exp}  Let $\Sigma$ be a $C^2$ hypersurface in a
vertically rigid sub-Riemannian manifold $M$.  Then, the {\bf horizontal
exponential surface} at $x \in \Sigma$ is defined to be the union of
all the horizontal curves in $\Sigma$ passing through $x$.  We denote
this subset of $\Sigma$ by $\Sigma_0(x)$.  
\enD

\bgD{Horizontal Tangent plane}   Let $\Sigma$ be a $C^2$ hypersurface in a
vertically rigid sub-Riemannian manifold $M$.  Then, the {\bf horizontal
tangent plane} at a noncharacteristic point $x \in \Sigma$, is defined as
\[T_x^h\Sigma = \{exp_x(v) | g(v,e_0(x)),v \in T_xM\}\]
where $exp$ is the Riemannian exponential map.
\enD
This horizontal tangent plane in a Carnot group can also be defined by blowing up the
metric at a given point (see \cite{FSSC2}).  

This gives us a geometric interpretation of these curvature conditions
analogous to the Riemannian setting:

\bgT{II-tangent plane}  Let $\Sigma$ be a $C^2$ hypersurface in $M$, a
Carnot group, and let $\{\kappa_i\}$ be the
set of horizontal principle curvatures of $\Sigma$ at a noncharacteristic point $x$.
Then, $\Sigma_0(x)$ locally lies to one side of $T_x^h\Sigma$ if and
only if the surface is horizontally positively curved at $x$.
Similarly, if $\Sigma$ is horizontally negatively curved at $x$, then any
neighborhood of $x$ in $\Sigma_0(x)$ intersects $T_x^h\Sigma$ at points other than $x$
\enT

\pf
First assume that $\Sigma_0(x)$ locally lies to one side of the
horizontal tangent plane at $x$.  Then, as any curves in $\Sigma_0(x)$
must also lie to one side of the tangent plane, we have that $\langle
\nabla_{\dot{c_1}}\dot{c_1}, e_0\rangle $ and $ \langle
\nabla_{\dot{c_2}}\dot{c_2}, e_0\rangle$ are either both non-positive or both non-negative at $x$ for any $c_1,c_2 \in
\Sigma_0(x)$.  Thus, $\Sigma$ is horizontally non-negatively curved at
$x$.  Conversely, if every curve $ c \in \Sigma_0(x)$ has positive
horizontal curvature, then $\langle \nabla_{\dot{c}}\dot(c),e_0\rangle$ is
either non-positive or non-negative.  Without loss of generality, we will
assume it to be non-negative.  But, geometrically, this says that with
respect to the connection $\nabla$, each curve $c$ locally changes in
the direction of $e_0$ and cannot move towards $-e_0$.  To finish the proof, consider a curve $c
\in \Sigma_0(x)$ in a small neighborhood of $x$.  If we let $v_i$ be the
left invariant unit vector field on $M$ so that $v_i(x)=e_i(x)$, we have that
$T^h_x\Sigma$ is the integral submanifold of the distribution
perpendicular to $v_0$ (with respect to the Riemannian metric).  With
respect to these vector fields, we write 
\[\dot{c}(t)= c_0(t)\; v_0 + \dots + c_N(t) \;v_N\]
where $c_0(0)=0$.  Then, computing with respect to the Riemannian
metric, we have
\begin{equation*}
\begin{split}
\langle \nabla_{\dot{c}} \dot{c},v_0\rangle &= \ddot{c}_0(t) +
\sum_{i,j} \dot{c}_i(t)\dot{c}_j(t) \langle \nabla_{v_i}v_j,v_0
\rangle \\
&= \ddot{c}_0(t)
\end{split}
\end{equation*}
The last equation follows because, as a Carnot group is graded, for left
invariant horizontal vector fields $v_1,v_2,v_3$, $\langle
[v_1,v_2],v_3 \rangle=0$.  Noticing that at $t=0$, we have 
\[\langle \nabla_{\dot{c}} \dot{c}(0),e_0\rangle =  \langle
\nabla_{\dot{c}} \dot{c}(0),v_0(c(0))\rangle = \ddot{c}_0(0)\]
The hypothesis that $II_0$ is positive semi-definite shows that
$\ddot{c}_0(0) \geq 0$ and the result follows.  A similar argument shows
the last statement.
\epf

This leads us to define a notion of convexity in sub-Riemannian
manifolds. 

\bgD{horiz convex}  A subset $U$ of a Carnot group $M$ with $C^2$ boundary $\Sigma$ is {\bf horizontally geometrically
  convex} (or hg-convex) if, at each noncharacteristic point $x \in \Sigma$,
$\Sigma_0(x)$ lies to one side of $T_x^h\Sigma$.  We say that $\Sigma$ is locally
hg-convex at a noncharacteristic point $x$ if there exists an $\epsilon >0$ so that $\Sigma_0(x)
\cap B(x,\epsilon)$ lies to one side of $T_x^h\Sigma$.
\enD

With this definition, we have yet another analogue of Euclidean
minimal surface theory:

\bgC{notconvex}  If $\Sigma$ is a $C^2$ minimal hypersurface in a
Carnot group then $\Sigma$ cannot bound a
locally hg-convex set.  
\enC

\pf
As the distribution $V_0$ is non-integrable, every $C^2$ surface must have at least one noncharacteristic point. The result then follows from \rfT{II-tangent plane} and \rfC{Negative Curvature}. 

\epf

 \rfT[intro]{Main 3} in the introduction is the combination of this corollary and \rfC{Negative Curvature}.

We note that this notion of convexity is distinct from the notion of set convexity
in Carnot groups introduced by Danielli, Garofalo and Nhieu in
\cite{DGN:convex} which is based on weakly convex defining functions
(the reader should also see the work of Lu, Manfredi and Straffolini \cite{LMS}
which independently presented a notion of convex functions at the same
time).  In particular, we point out that there are minimal surfaces in
the first Heisenberg group which are the boundary of a convex region
in the sense of \cite{DGN:convex} (for example, the plane $t=0$ and the surface
$t=\frac{xy}{2}$) and others that bound nonconvex sets (for example,
$t= \frac{xy}{2}-\frac{x^3}{3}$).  Moreover, we further point out that
as $\phi(x,y,t)= \frac{xy}{2}+g(x)-t$ has minimal level sets in
$\mathbb{H}$ for every
choice of $C^2$ function $g$, we can easily produce many functions
that violate the various convexity conditions presented in
\cite{DGN:convex} and \cite{LMS}.  As the level sets of $\phi=x$ are minimal in $\mathbb{H}$ and satisfy all of the convexity
properties defined in \cite{DGN:convex}, we see that none of these
notions can be made equivalent to the notion of hg-convexity.

\bgR{Cake}  We remark that the definitions and proofs in this section
are direct generalizations or adaptations of the Euclidean and/or
Riemannian machinery.  While the proofs are quite straightforward, we
point out that this is due to a correct choice of geometric structure,
in this case the adapted connection, that allows for the ease of the
proofs.  Without such machinery, the statements about the relation
between horizontal curvature and minimality/isoperimetry above were known only in certain
three dimensional sub-Riemannian manifolds where they reduce to much
simpler statements.  
\enR